\documentclass[titlepage,12pt]{article}

\usepackage[pctex32]{graphics}

\def\be{\begin{equation}}

\def\ee{\end{equation}}

\usepackage{amssymb}

\begin{document}

\begin{center}
{\Large{\bf Orthogonal Asymptotic Lines \\ on Surfaces Immersed in
$\mathbf R ^4$}}
\end{center}

\vspace{2cm}

\begin{center}
{\large Luis Fernando Mello}
\end{center}

\begin{center}
{\em Instituto de Ci\^encias Exatas \\Universidade Federal de
Itajub\'a \\CEP 37.500-903, Itajub\'a, MG, Brazil
\\}e--mail:lfmelo@unifei.edu.br
\end{center}

\vspace{3cm}

\begin{center}
{\bf Abstract}
\end{center}

\vspace{0.5cm}

In this paper we study some properties of surfaces immersed in
$\mathbf R ^4$ whose asymptotic lines are orthogonal. We also
analyze necessary and sufficient conditions for the
hypersphericity of surfaces in $\mathbf R ^4$.

\vspace{3cm}

\noindent {\bf Key words and phrases}: asymptotic lines,
inflection points, normal curvature, hypersphericity.

\vspace{0.3cm}

\noindent 2000 Mathematics Subject Classification: 53A05, 58C25.

\vspace{0.3cm}

\noindent Research partially supported by FAPEMIG grant EDT
1929/03.

\newpage
\baselineskip=20pt

\section{Introduction}\label{S:1}

There are two different ways to construct line fields on surfaces
immersed in $\mathbf R ^4$. The first one consists in considering
the ellipse of curvature in the normal bundle of the surface and
taking the pull back of points on this ellipse to define tangent
direction fields. Examples of this approach are given by: the {\it
lines of axial curvature}, along which the second fundamental form
points in the direction of the large and the small axes of the
ellipse of curvature; the {\it mean directionally curved lines},
along which the second fundamental form points in the direction of
the mean curvature vector; and the {\it asymptotic lines}, along
which the second fundamental form points in the direction of the
tangent lines to the ellipse of curvature.

The other way consists in defining the {\it $\nu$-principal
curvature lines}, along which the surface bends extremally in the
direction of the normal vector $\nu$. To this end, we need to take
an unitary normal vector field $\nu$ and follow the classical
approach for surfaces immersed in $\mathbf R ^3$.

The lines of axial curvature are globally defined and their
singularities are the axiumbilic points where the ellipse of
curvature becomes either a circle or a point. The axiumbilic
points and the lines of axial curvature are assembled into two
axial configurations. The first one is defined by the axiumbilics
and the field of orthogonal lines on which the surface is curved
along the large axis of the ellipse of curvature. The second one
is defined by the axiumbilics and the field of orthogonal lines on
which the surface is curved along the small axis of the ellipse of
curvature. Each axial configuration is a net consisting of
orthogonal curves and axiumbilic points. Therefore a line of axial
curvature is not necessarily a simple regular curve; it can be
immersed with transversal crossings. The differential equation of
lines of axial curvature is a quartic differential equation
according to \cite{GS, GGTG1, GGTG2}. A global analysis of the
lines of axial curvature was developed in \cite{GS}.

The mean directionally curved lines are globally defined and their
singularities are either the inflection points, where the ellipse
of curvature is a radial line segment, or the minimal points,
where the mean curvature vector vanishes. It was shown in \cite{M}
that the differential equation of mean directionally curved lines
fits into the class of quadratic or binary differential equations.
The global behavior of mean directionally curved lines was studied
in \cite{M}.

The asymptotic lines do not need to be globally defined on the
surfaces and in general are not orthogonal. It was shown in
\cite{MRR} that a necessary and sufficient condition for existence
of the globally defined asymptotic lines on a surface
$\mathbf{M}^2$ in $\mathbf R ^4$ is the local convexity of
$\mathbf{M}^2$. The differential equation of asymptotic lines is
also a quadratic differential equation and their singularities are
the inflection points.

The $\nu$-principal curvature lines are orthogonal and globally
defined on surfaces immersed in $\mathbf R ^4$ and their
singularities are the $\nu$-umbilic points, where the
$\nu$-principal curvatures coincide. The differential equation of
$\nu$-principal curvature lines is a quadratic differential
equation according to \cite{SR}. An analysis of $\nu$-principal
curvature lines near generic $\nu$-umbilic points is presented in
\cite{SR} and in \cite{GSa} the $\nu$-principal cycles (closed
$\nu$-principal curvature lines) are studied. A global analysis of
the $\nu$-principal curvature lines was developed in \cite{GMS},
for $\nu = H$, where $H$ is the normal mean curvature vector.

We prove in \cite{M2} that the orthogonality of the asymptotic
lines is equivalent to the vanishing of the normal curvature. This
result has been already obtained by Romero-Fuster and
S\'anchez-Bringas in \cite{RS} using a different approach. We also
prove in \cite{M2} that the quartic differential equation of lines
of axial curvature can be written as the product of the quadratic
differential equations of mean directionally curved lines and
asymptotic lines if and only if the normal curvature of $\alpha$
vanishes at every point. Thus if the normal curvature of $\alpha$
vanishes at every point then the axial curvature cross fields
split into four direction fields and therefore it is not possible
that the lines of axial curvature have transversal crossings.

On the other hand, it is well known that a point $p$ is
semiumbilic if and only if the normal curvature vanishes at $p$,
\cite{RS}. Semiumbilic points are interesting from the viewpoint
of the theory of singularities of functions. Observe now that we
have analogous statements if instead of vanishing normal curvature
it is required semiumbilicity.

We say that an immersion $\alpha:\mathbf{M}^2\rightarrow \mathbf R
^4$ is {\it hyperspherical} if its image is contained in a
hypersphere. In this work we study some properties of surfaces
immersed in $\mathbf R ^4$ whose asymptotic lines are orthogonal.
In particular, we relate the property of having globally defined
orthogonal asymptotic lines with hypersphericity, obtaining the
following theorem.

\noindent {\bf Theorem 3.2.} Let $\alpha:\mathbf{M}^2\rightarrow
\mathbf R ^4$ be an immersion of a smooth oriented surface with
globally defined orthogonal asymptotic lines. Suppose that there
exist an unitary normal vector field $\nu$ and $r > 0$ such that
the distance from the projection of the ellipse of curvature
$\varepsilon_{\alpha}(p)$ onto the $\nu$-axis to $p$ is $r$, for
all $p \in \mathbf{M}^2$, and the Gaussian curvature $K \neq r^2$.
Then $\alpha$ is hyperspherical.

Finally, theorem 3.4 of \cite{RS}, lemma 2.1 and theorem 2.1 of
\cite{M2} and results of this paper are put together in Theorem
\ref{teo5} establishing seven other equivalent conditions to the
orthogonality of the asymptotic lines.

This paper is organized as follows. A review of properties of the
first and second fundamental forms, the ellipse of curvature and
the line fields on surfaces immersed in $\mathbf R ^4$ is
presented in section \ref{S:2}. General aspects of the curvature
theory for surfaces immersed in $\mathbf R ^4$ are presented in
the works of Forsyth \cite{F}, Wong \cite{W}, Little \cite{L} and
Asperti \cite{A}. Section \ref{S:3} is devoted to the study of
orthogonal asymptotic lines as well as hypersphericity of
immersions. Finally, in section \ref{S:4} some general problems
are stated.

\newtheorem{teo}{Theorem}[section]
\newtheorem{lema}[teo]{Lemma}
\newtheorem{prop}[teo]{Proposition}
\newtheorem{cor}[teo]{Corollary}

\section{%2.
Line fields on surfaces in $\mathbf R ^4$}\label{S:2} For sake of
completeness in this section we present a survey of the relevant
notions that will need later. Let $\alpha:\mathbf{M}^2\rightarrow
\mathbf R ^4$ be an immersion of a smooth oriented surface into
$\mathbf R ^4$, which is endowed with the Euclidean inner product
$\langle \cdot,\cdot \rangle$ and is oriented. In this paper
immersions are assumed to be $C^\infty$. Denote respectively by
$\mathbf{TM}$ and $\mathbf{NM}$ the tangent and the normal bundles
of $\alpha$ and by $T_{p}\mathbf{M}$ and $N_{p}\mathbf{M}$ the
respective fibers, i.e., the tangent and the normal planes at
$p\in \mathbf{M}^2$. Let $\{\nu_1,\nu_2\}$ be a frame of vector
fields orthonormal to $\alpha$. Assume that $(u,v)$ is a positive
chart of $\mathbf{M}^2$ and that
$\{\alpha_u,\alpha_v,\nu_1,\nu_2\}$ is a positive frame of
$\mathbf R ^4$. In such a chart $(u,v)$ the first fundamental form
of $\alpha$, $I_\alpha$, is given by
\[
I=I_\alpha=\langle d\alpha,d\alpha \rangle =E du^2+2F dudv+G dv^2,
\]
where $E=\langle \alpha_u,\alpha_u \rangle$, $F=\langle
\alpha_u,\alpha_v \rangle$ and $G=\langle \alpha_v,\alpha_v
\rangle$. The second fundamental form of $\alpha$, $II_\alpha$, is
defined in terms of the $\mathbf{NM}$-valued quadratic form
\[
II=II_\alpha=\langle d^2\alpha,\nu_1 \rangle \nu_1+\langle
d^2\alpha,\nu_2 \rangle \nu_2=II_{\nu_1}\nu_1+II_{\nu_2}\nu_2,
\]
where
\[
II_{\nu_i} = II_{\nu_i,\alpha}=e_{i}du^2+2f_{i}dudv+g_{i}dv^2,
\]
$e_{i}=\langle \alpha_{uu},\nu_{i} \rangle$, $f_{i}=\langle
\alpha_{uv},\nu_{i} \rangle$, and $g_{i}=\langle
\alpha_{vv},\nu_{i} \rangle$, for $i=1,2$.

The following functions are associated to $\alpha$ (see \cite{L}):
\begin{enumerate}

\item The {\it mean curvature vector} of $\alpha$
\[
H=H_{\alpha}=H_{1}\nu_{1}+H_{2}\nu_{2},
\]
where
\[
H_{i}=H_{i,\alpha}=\frac{E g_{i}-2 F f_{i}+G e_{i}}{2(E G-F^2)},
\]
for $i=1,2$;

\item The {\it normal curvature} of $\alpha$
\[
k_{N}=k_{N,\alpha}=\frac{E(f_{1}g_{2}-f_{2}g_{1})-F(e_{1}g_{2}-e_{2}g_{1})+G(e_{1}f_{2}-e_{2}f_{1})
}{2(E G-F^2)};
\]

\item The {\it resultant} ${\it \Delta}$ of $II_{1,\alpha}$ and
$II_{2,\alpha}$
\[
{\it \Delta}={\it \Delta}_\alpha={1\over {4(E G-F^2)}}
\left|\begin{array}{cccc}
                                    e_{1} & 2f_{1} & g_{1} & 0 \\
                                    e_{2} & 2f_{2} & g_{2} & 0 \\
                                    0     & e_{1}  & 2f_{1}& g_{1} \\
                                    0     & e_{2}  & 2f_{2}& g_{2}
\end{array} \right|;
\]

\item The {\it Gaussian curvature} of $\alpha$
\[
K=K_{\alpha}=\frac{e_1 g_1 - (f_1)^2 + e_2 g_2 - (f_2)^2}{E
G-F^2};
\]

\item The {\it normal curvature vector} of $\alpha$ defined by
$\eta(p,v)=\frac{II(p,v)}{I(p,v)}$.
\end{enumerate}

The image of the unitary tangent circle $\mathbf{S}^1$ by
$\eta(p):T_p\mathbf{M}\rightarrow N_p\mathbf{M}$ describes an
ellipse in $N_p\mathbf{M}$ called {\it ellipse of curvature} of
$\alpha$ at {\it p} and denoted by $\varepsilon_{\alpha}(p)$. This
ellipse may degenerate into a line segment, a circle or a point.
The center of the ellipse of curvature is the mean curvature
vector {\it H} and the area of $\varepsilon_{\alpha}(p)$ is given
by ${\pi\over 2} \left|k_N(p) \right|$. The map $\eta(p)$
restricted to $\mathbf{S}^1$, being quadratic, is a double
covering of the ellipse of curvature. Thus every point of the
ellipse corresponds to two diametrically opposed points of the
unitary tangent circle. The ellipse of curvature is invariant by
rotations in both the tangent and normal planes.

A point $p\in \mathbf{M}^2$ is called a {\it minimal point} of
$\alpha$ if $H(p)=0$ and it is called an {\it inflection point} of
$\alpha$ if ${\it \Delta}(p)=0 $ and $k_N(p)=0$. It follows that
$p\in \mathbf{M}^2$ is an inflection point if and only if its
ellipse of
curvature is a radial line segment \cite{L}.\\
\noindent{\bf Lines of axial curvature.} The four vertices of the
ellipse of curvature $\varepsilon_{\alpha}(p)$ determine eight
points on the unitary tangent circle which define two crosses in
the tangent plane. Thus we have two cross fields on $\mathbf{M}^2$
called {\it axial curvature cross fields}. This construction fails
at the {\it axiumbilic points} where the ellipse of curvature
becomes either a circle or a point. Generically the index of an
isolated axiumbilic point is $\pm {1\over 4}$ (see \cite{GS,
GGTG1, GGTG2}). The integral curves of the axial curvature cross
fields are the {\it lines of axial curvature}.

Generically there is no good way to distinguish one end of the
large (or small) axis of $\varepsilon_{\alpha}(p)$ and therefore
pick out a direction of the cross field. Thus a line of axial
curvature is not necessarily a simple regular curve; it can be
immersed with transversal crossings.

The differential equation of the lines of axial curvature is a
quartic differential equation of the form
\begin{equation} \label{1}
Jac\biggl(\|\eta-H\|^2,I \biggr)=0,
\end{equation}
\noindent where
\[
Jac(\cdot, \cdot)={{\partial(\cdot, \cdot)}\over {\partial
(du,dv)}},
\]
which according to \cite{GS} can be written as
\begin{equation} \label{2}
A_0du^4+A_1du^3dv+A_2du^2dv^2+A_3dudv^3+A_4dv^4=0,
\end{equation}
\noindent where
\[
A_0=a_0E^3, \; A_1=a_1E^3, \; A_2=-6a_0GE^2+3a_1FE^2,
\]
\[
A_3=-8a_0EFG+a_1E(4F^2-EG), \; A_4=a_0G(EG-4F^2)+a_1F(2F^2-EG),
\]
\[
a_0=4\biggl[F(EG-2F^2)(e_1^2+e_2^2)-Ea_6a_2-E^2F(a_3+a_5)+E^3a_4\biggr],
\]
\[
a_1=4\biggl[Ga_6(e_1^2+e_2^2)+8EFGa_2+E^3(g_1^2+g_2^2)-2E^2G(a_3+a_5)\biggr],
\]
\[
a_2=e_1f_1+e_2f_2, \; a_3=e_1g_1+e_2g_2, \; a_4=f_1g_1+f_2g_2,
\]
\[
a_5=2(f_1^2+f_2^2), \; a_6=EG-4F^2.
\]
\noindent{\bf Mean directionally curved lines.} The line through
the mean curvature vector $H(p)$ meets $\varepsilon_{\alpha}(p)$
at two diametrically opposed points. This construction induces two
orthogonal directions on $T_p\mathbf{M}^2$. Therefore we have two
orthogonal direction fields on $\mathbf{M}^2$ called {\it
H-direction fields}. The singularities of these fields, called
here {\it H-singularities}, are the points where either $H=0$
(minimal points) or at which the ellipse of curvature becomes a
radial line segment (inflection points). Generically the index of
an isolated {\it H}-singularity is $\pm {1\over 2}$ \cite{M}. The
integral curves of the {\it H}-direction fields are the {\it mean
directionally curved lines}.

The differential equation of mean directionally curved lines is a
quadratic differential equation of the form \cite{M}
\begin{equation} \label{3}
Jac\{Jac(II_{\nu_1},II_{\nu_2}),I\}=0,
\end{equation}
\noindent which can be written as
\begin{equation} \label{4}
B_1(u,v)du^2+2B_2(u,v)dudv+B_3(u,v)dv^2=0,
\end{equation}
\noindent where
\[
B_1=(e_1g_2-e_2g_1)E+2(e_2f_1-e_1f_2)F,\:\:
B_2=(f_1g_2-f_2g_1)E+(e_2f_1-e_1f_2)G,
\]
\[
B_3=2(f_1g_2-f_2g_1)F+(e_2g_1-e_1g_2)G.
\]
\noindent{\bf Asymptotic lines.} Suppose that $p$ (the origin of
$N_p\mathbf{M}^2$) lies outside $\varepsilon_{\alpha}(p)$, for all
$p\in \mathbf{M}^2$. The two points on $\varepsilon_{\alpha}(p)$
at which the lines through the normal curvature vectors are
tangent to $\varepsilon_{\alpha}(p)$ induce a pair of directions
in $T_p\mathbf{M}^2$ which in general are not orthogonal. Thus we
have two tangent direction fields on $\mathbf{M}^2$, called {\it
asymptotic direction fields}. The singularities of these fields
are the points where the ellipse of curvature becomes a radial
line segment, i.e., the {\it inflection points}. Generically the
index of an isolated inflection point is $\pm {1\over 2}$
\cite{GMRR}. The integral curves of the asymptotic direction
fields are the {\it asymptotic lines}.

The differential equation of asymptotic lines is a quadratic
differential equation of the form \cite{M}
\begin{equation} \label{5}
Jac(II_{\nu_1},II_{\nu_2})=0,
\end{equation}
\noindent which can be written as
\begin{equation} \label{6}
T_1(u,v)du^2+T_2(u,v)dudv+T_3(u,v)dv^2=0,
\end{equation}
\noindent where
\[
T_1=e_1f_2-e_2f_1, \; T_2=e_1g_2-e_2g_1, \; T_3=f_1g_2-f_2g_1.
\]
\noindent {\bf $\nu$-Principal curvature lines.} The projection of
the pullback, ${\alpha}{^*} (\mathbf R ^4)$, of the tangent bundle
of $\mathbf R ^4$ onto the tangent bundle of an immersion $\alpha$
will be denoted by $\Pi_{\alpha,T}$. This vector bundle is endowed
with the standard metric induced by the Euclidean one in $\mathbf
R ^4$.

Denote by $\nu = \nu_{\alpha}$ the {\it unit normal vector field}
of ${\alpha}$. The  eigenvalues $k_1 = k_{1,\alpha} \leq
k_{2,\alpha} = k_2$ of the {\it Weingarten operator} ${\mathcal
W}_\alpha = -\Pi_{\alpha,T}D\nu_\alpha $ of $\mathbf {TM}$ are
called the {\it $\nu$-principal curvatures} of $\alpha$. The
points where $k = k_1 = k_2$ will be called the {\it
$\nu$-umbilic} points of $\alpha$ and define the set ${\mathcal
S}_u = {\mathcal S}_{u,\alpha}$. We say that $\alpha$ is {\it
$\nu$-umbilical} if each point of the immersion is $\nu$-umbilic.
Outside ${\mathcal S}_u$ are defined the {\it minimal},
$L_{m,\alpha}$, and the {\it maximal}, $L_{M,\alpha}$, {\it
$\nu$-principal line fields} of $\alpha$, which are the
eigenspaces of ${\mathcal W}_\alpha$ associated respectively to
$k_1$ and $k_2$. Generically the index of an isolated
$\nu$-umbilic point is $\pm {1\over 2}$ \cite{SR}. The integral
curves of the $\nu$-principal line fields are the {\it
$\nu$-principal curvature lines}.

In a local chart $(u,v)$ the $\nu$-principal curvatures lines are
characterized as the solutions of the following quadratic
differential equation \cite{SR}
\begin{equation} \label{pcl}
(Fg_\nu - f_\nu G)dv^2+(E g_\nu - e_\nu G)dudv+(E f_\nu - F
e_\nu)du^2=0,
\end{equation}
\noindent where $E$, $F$ and $G$ are the coefficients of the first
fundamental form and  $e_\nu = <\alpha_{uu}, \nu>$, $f_\nu =
<\alpha_{uv},\nu>$ and $g_\nu = <\alpha_{vv},\nu>$ are the
coefficients of the {\it second fundamental form relative to}
$\nu$, denoted by $II_\nu = II_{\nu_\alpha}$. Equation (\ref{pcl})
is equivalently written as
\begin{equation} \label{pcl2}
Jac (II_\nu, I) = 0.
\end{equation}

\section{%3.
Orthogonal asymptotic lines}\label{S:3} Let
$\alpha:\mathbf{M}^2\rightarrow \mathbf R ^4$ be an immersion of a
smooth oriented surface into $\mathbf R ^4$. In \cite{GS} Garcia
and Sotomayor prove the following theorem: Suppose that the image
of the surface $\mathbf{M}^2$ by $\alpha$ is contained into
$\mathbf R ^3$. Then the quartic differential equation of lines of
axial curvature is the product of the quadratic differential
equation of its principal curvature lines and the quadratic
differential equation of its mean curvature lines. It is
interesting to observe that every point of $\mathbf{M}^2$ is an
inflection point.

We have established in \cite{M} the following theorem: Let
$\alpha:\mathbf{M}^2\rightarrow \mathbf{S}^3(r)$ be an immersion
of a smooth oriented surface into a 3-dimensional sphere of radius
$r>0$. Consider the natural inclusion
$i:\mathbf{S}^3(r)\rightarrow \mathbf R ^4$ and the composition
$i\circ \alpha$ also denoted by $\alpha$. Then the quartic
differential equation of lines of axial curvature (\ref{1}) can be
written as
\begin{equation} \label{7}
Jac\{Jac(II_{\nu_1},II_{\nu_2}),I\} \cdot
Jac(II_{\nu_1},II_{\nu_2})=0,
\end{equation}
\noindent where the first expression in (\ref{7}) is the quadratic
differential equation of mean directionally curved lines (\ref{3})
and the second one is the quadratic differential equation of
asymptotic lines (\ref{5}).

It is interesting to observe that in the above construction the
asymptotic lines are orthogonal and the normal curvature of
$\alpha$ vanishes at every point. This is a particular case of the
following theorem proved in \cite{M2}, which was also obtained in
\cite{RS} using a different approach: Let
$\alpha:\mathbf{M}^2\rightarrow \mathbf R ^4$ be an immersion of a
smooth oriented surface with isolated inflection points. The
immersion $\alpha$ has orthogonal asymptotic lines if and only if
the normal curvature of $\alpha$ vanishes at every point.

We have established in \cite{M2} the following theorem: Let
$\alpha:\mathbf{M}^2\rightarrow \mathbf R ^4$ be an immersion of a
smooth oriented surface with isolated inflection points. The
quartic differential equation of lines of axial curvature
(\ref{1}) can be written as
\begin{equation} \label{8}
Jac\{Jac(II_{\nu_1},II_{\nu_2}),I\}\cdot
Jac(II_{\nu_1},II_{\nu_2})=0,
\end{equation}
\noindent where the first expression in (\ref{8}) is the quadratic
differential equation of mean directionally curved lines (\ref{3})
and the second one is the quadratic differential equation of
asymptotic lines (\ref{5}), if and only if the normal curvature of
$\alpha$ vanishes at every point.

We can prove the following corollary: Let
$\alpha:\mathbf{M}^2\rightarrow \mathbf R ^4$ be an immersion of a
smooth oriented surface into $\mathbf R ^4$. If the immersion
$\alpha$ has orthogonal asymptotic lines then the inflection
points are obtained where the ellipse of curvature becomes a
point. In fact, from Equation (\ref{8})
\begin{equation} \label{9}
Jac\biggl(\|\eta-H\|^2,I
\biggr)=Jac\{Jac(II_{\nu_1},II_{\nu_2}),I\}\cdot
Jac(II_{\nu_1},II_{\nu_2})=0.
\end{equation}
\noindent As the inflection points are singularities of asymptotic
lines then by (\ref{9}) they are singularities of lines of axial
curvature. But the singularities of lines of axial curvature are
the points where the ellipse of curvature becomes either a circle
or a point. Thus the only possibility in this case is that the
ellipse of curvature becomes a point.

\begin{teo} \label{teo1}
Let $\alpha:\mathbf{M}^2\rightarrow \mathbf{S}^3(r)$ be an
immersion of a smooth oriented surface into a 3-dimensional sphere
of radius $r>0$. Consider the natural inclusion
$i:\mathbf{S}^3(r)\rightarrow \mathbf R ^4$ and the composition
$i\circ \alpha$ also denoted by $\alpha$. Then there exist an
unitary normal vector field $\nu$ and $\lambda > 0$ such that the
ellipse of curvature $\varepsilon_{\alpha}(p)$ is a line segment
with the following property: the distance from the projection of
$\varepsilon_{\alpha}(p)$ onto the $\nu$-axis to $p$ is $\lambda$,
for all $p \in \mathbf{M}^2$.
\end{teo}

\noindent{\bf Proof.} Let $\{\nu_1,\nu_2\}$ be a frame of vector
fields orthonormal to $\alpha$, where $\nu_1(p)\in
T_{p}\mathbf{S}^3(r)$ and $\nu_2(p)$ is the inward normal to
$\mathbf{S}^3(r)$, for all $p\in \mathbf{M}^2$. Thus
\[
\nu_2\equiv - {\frac{1} {r} \; {\alpha}},\; e_2 = \frac{1} {r} \;
E,\; f_2 = \frac{1} {r}\; F\; \mbox{and}\; g_2 = \frac{1} {r} \;
G,
\]
where $E$, $F$ and $G$ are the coefficients of the first
fundamental form of $\alpha$. It follows that
\[
II_{\nu_2} = \frac{1} {r} \; I.
\]
Now
\[
\eta=\frac{II}{I}=\frac{II_{\nu_1}}{I} \; \nu_1 +
\frac{II_{\nu_2}}{I} \; \nu_2= \frac{II_{\nu_1}}{I} \; \nu_1 +
\frac{1} {r} \; \nu_2.
\]
This implies that the ellipse of curvature
$\varepsilon_{\alpha}(p)$ is a line segment orthogonal to $\nu_2$,
for all $p\in \mathbf{M}^2$. Define $\nu = \nu_2$ and $\lambda =
\frac{1}{r}$. The theorem is proved.

Let $\alpha:\mathbf{M}^2\rightarrow \mathbf R ^4$ be an immersion
of a smooth oriented surface with globally defined orthogonal
asymptotic lines. Then the ellipse of curvature
$\varepsilon_{\alpha}(p)$ is a line segment for all $p \in
\mathbf{M}^2$ except at the inflection points. We say that the
immersion $\alpha$ has {\it constant projection} if there exist an
unitary normal vector field $\nu$ and $r > 0$ such that the
distance from the projection of $\varepsilon_{\alpha}(p)$ onto the
$\nu$-axis to $p$ (the origin of $N_p \mathbf{M}^2$) is $r$, for
all $p \in \mathbf{M}^2$. The constant $r$ is called {\it distance
of projection}.

Theorem \ref{teo1} shows that if $\alpha$ is hyperspherical then
$\alpha$ has constant projection whose distance of projection is
$r^{-1}$, where $r$ is the radius of the hypersphere. The converse
of Theorem~\ref{teo1} is given by the following theorem.

\begin{teo} \label{teo2}
Let $\alpha:\mathbf{M}^2\rightarrow \mathbf R ^4$ be an immersion
of a smooth oriented surface with globally defined orthogonal
asymptotic lines. Suppose that $\alpha$ has constant projection
with distance of projection $r > 0$, and the Gaussian curvature $K
\neq r^2$. Then $\alpha$ is hyperspherical.
\end{teo}

\noindent{\bf Proof.} Since all the notions of this paper are
independents of the chart it is enough to prove this theorem for
an orthogonal one. By hypothesis there is an unitary normal vector
field $\nu$ orthogonal to $\varepsilon_{\alpha}(p)$, for all $p
\in {\mathbf{M}^2}$. We can take $\{\nu_1 = \nu^\perp,\nu_2 = \nu
\}$ a frame of vector fields orthonormal to $\alpha$, where
$\nu^\perp(p)$ is parallel to $\varepsilon_{\alpha}(p)$, such that
$\{\alpha_u,\alpha_v,\nu^\perp,\nu\}$ is a positive frame of
$\mathbf R ^4$, for a positive orthogonal chart $(u,v)$ of
$\mathbf{M}^2$. Thus $e_2 = r E, \; f_2 = 0, \; g_2 = r G$. The
immersion $\alpha$ satisfies the Codazzi equations \cite{F}
\begin{equation}\label{c1}
(e_1)_v-(f_1)_u = \Gamma_{12}^1 e_1 + \left( \Gamma_{12}^2 - \Gamma_{11}^1 \right)f_1-
\Gamma_{11}^2 g_1 - a_{12}^3 e_2 + a_{11}^3 f_2,
\end{equation}
\begin{equation}\label{c2}
(e_2)_v-(f_2)_u = \Gamma_{12}^1 e_2 + \left( \Gamma_{12}^2 - \Gamma_{11}^1 \right)f_2-
\Gamma_{11}^2 g_2 - a_{12}^3 e_1 + a_{11}^3 f_1,
\end{equation}
\begin{equation}\label{c3}
(f_1)_v-(g_1)_u = \Gamma_{22}^1 e_1 + \left( \Gamma_{22}^2 - \Gamma_{12}^1 \right)f_1-
\Gamma_{12}^2 g_1 + a_{12}^3 f_2 - a_{11}^3 g_2,
\end{equation}
\begin{equation}\label{c4}
(f_2)_v-(g_2)_u = \Gamma_{22}^1 e_2 + \left( \Gamma_{22}^2 - \Gamma_{12}^1 \right)f_2-
\Gamma_{12}^2 g_2 - a_{12}^3 f_1 + a_{11}^3 g_1,
\end{equation}
\noindent and the following structure equations \cite{F}
\begin{equation}\label{s1}
(\nu^\perp)_u = a_{11}^1 \alpha_u + a_{11}^2 \alpha_v + a_{11}^3
\nu,
\end{equation}
\begin{equation}\label{s2}
(\nu^\perp)_v = a_{12}^1 \alpha_u + a_{12}^2 \alpha_v + a_{12}^3
\nu,
\end{equation}
\begin{equation}\label{s3}
\nu_u = a_{21}^1 \alpha_u + a_{21}^2 \alpha_v - a_{11}^3 \nu^\perp,
\end{equation}
\begin{equation}\label{s4}
\nu_v = a_{22}^1 \alpha_u + a_{22}^2 \alpha_v - a_{12}^3 \nu^\perp,
\end{equation}
\noindent where
\[
a_{11}^1 = \frac{f_1F-e_1G}{EG-F^2}, \; a_{11}^2 = \frac{e_1F-f_1E}{EG-F^2}, \;
a_{12}^1 = \frac{g_1F-f_1G}{EG-F^2}, \; a_{12}^2 = \frac{f_1F-g_1E}{EG-F^2},
\]
\[
a_{21}^1 = \frac{f_2F-e_2G}{EG-F^2}, \; a_{21}^2 = \frac{e_2F-f_2E}{EG-F^2}, \;
a_{22}^1 = \frac{g_2F-f_2G}{EG-F^2}, \; a_{22}^2 = \frac{f_2F-g_2E}{EG-F^2},
\]
\noindent and $\Gamma_{ij}^k$ are the Christoffel symbols of $\alpha$ \cite{F},
$i,j,k = 1, 2$, which in this case are given by
\[
\Gamma_{11}^1=\frac {E_u}{2E}, \; \Gamma_{11}^2= - \frac
{E_v}{2G}, \; \Gamma_{12}^1=\frac {E_v}{2E}, \;
\Gamma_{12}^2=\frac {G_u}{2G}, \; \Gamma_{22}^1= - \frac
{G_u}{2E}, \; \Gamma_{22}^2=\frac {G_v}{2G}.
\]
\noindent Substituting the above Christoffel symbols in the
Codazzi equations (\ref{c2}) and (\ref{c4}) we have respectively
\begin{equation}\label{c5}
rE_v = \frac{E_v}{2E}rE + \frac{E_v}{2G}rG - a_{12}^3 e_1 +
a_{11}^3 f_1
\end{equation}
\noindent and
\begin{equation}\label{c6}
-rG_u = - \frac{G_u}{2E}rE - \frac{G_u}{2G}rG - a_{12}^3 f_1 +
a_{11}^3 g_1.
\end{equation}
\noindent But Equations (\ref{c5}) and (\ref{c6}) are equivalents
to
\begin{equation}\label{c7}
-a_{12}^3 e_1 + a_{11}^3 f_1 = 0
\end{equation}
\noindent and
\begin{equation}\label{c8}
-a_{12}^3 f_1 + a_{11}^3 g_1 = 0,
\end{equation}
\noindent respectively. Now the Gaussian curvature is
\[
K = \frac{e_1 g_1 -(f_1)^2}{EG} + \frac{e_2 g_2}{EG} = \frac{e_1
g_1 -(f_1)^2}{EG} + r^2.
\]
By hypothesis $K  \neq r^2$, and thus
\begin{equation}\label{ga}
e_1 g_1 - (f_1)^2 \neq 0.
\end{equation}
\noindent From the Equations (\ref{c7}), (\ref{c8}) and (\ref{ga})
we have that
\begin{equation}\label{s5}
a_{11}^3=a_{12}^3 = 0.
\end{equation}
\noindent Substituting Equation (\ref{s5}) in (\ref{s3}) and
(\ref{s4}) results that
\[
\nu_u = -r \alpha_u \; \mbox{and}\; \nu_v = -r \alpha_v.
\]
Thus
\[
\nu = -r \alpha + \gamma,
\]
where $\gamma$ is a constant vector. Therefore
\[
\alpha = \frac {\gamma}{r} - \frac {1}{r} \nu.
\]
This means that $\alpha(\mathbf{M}^2)$ belongs to a hypersphere
with center $\frac{\gamma}{r}$ and radius $\frac{1}{r}$. The
theorem is proved.

The proof of the following theorem is immediate from the proof of
Theorem \ref{teo1}.

\begin{teo} \label{teo3}
Let $\alpha:\mathbf{M}^2\rightarrow \mathbf{S}^3(r)$ be an
immersion of a smooth oriented surface into a 3-dimensional sphere
of radius $r>0$. Consider the natural inclusion
$i:\mathbf{S}^3(r)\rightarrow \mathbf R ^4$ and the composition
$i\circ \alpha$ also denoted by $\alpha$. Then there exist an
unitary normal vector field $\nu$ and $\lambda > 0$ such that
$II_\nu = \langle d^2 \alpha,\nu \rangle = \lambda I $.
\end{teo}

The converse of Theorem \ref{teo3} is given by the following
theorem.

\begin{teo} \label{teo4}
Let $\alpha:\mathbf{M}^2\rightarrow \mathbf R ^4$ be an immersion
of a smooth oriented surface. Suppose that $\nu$ is an unitary
normal vector field such that $II_\nu = \langle d^2 \alpha,\nu
\rangle = \lambda I $, where $\lambda$ is a nonzero constant, and
the Gaussian curvature $K \neq \lambda^2$. Then $\alpha$ is
hyperspherical.
\end{teo}

\noindent {\bf Proof.} Take the positive frame
$\{\alpha_u,\alpha_v,\nu^\perp,\nu \}$. As $II_\nu = \langle d^2
\alpha,\nu \rangle = \lambda I $ we have
\[
\eta=\frac{II}{I}=\frac{II_{\nu^\perp}}{I} \; \nu^\perp +
\frac{II_\nu}{I} \; \nu = \frac{II_{\nu^\perp}}{I} \; \nu^\perp +
\lambda \; \nu.
\]
This implies that the ellipse of curvature
$\varepsilon_{\alpha}(p)$ is a line segment whose distance from
their projection onto the $\nu$-axis to $p$ is constant and equal
to $\lambda$, for all $p\in \mathbf{M}^2$. Therefore $\alpha$ has
constant projection with distance of projection $\lambda > 0$. As
$K \neq \lambda^2$ the theorem follows from Theorem \ref{teo2}.

Let $\alpha:\mathbf{M}^2\rightarrow \mathbf R ^4$ be an immersion
of a smooth oriented surface with globally defined orthogonal
asymptotic lines. Then the normal curvature of $\alpha$ vanishes
at every point. So there exist normal vector fields $\nu$ and
$\nu^\perp$ such that
\[
\eta=\frac{II}{I}=\frac{II_{\nu^\perp}}{I} \; \nu^\perp +
\frac{II_\nu}{I} \; \nu = \frac{II_{\nu^\perp}}{I} \; \nu^\perp +
\lambda \; \nu.
\]
\noindent Thus $II_\nu = \lambda I $, where $\lambda$ is a
positive scalar function on $\mathbf{M}^2$. This implies that
$\alpha$ is $\nu$-umbilical. The differential equation of
asymptotic lines (\ref{5}) is given by
\[
0 = Jac(II_{\nu^\perp},II_{\nu}) = Jac(II_{\nu^\perp},\lambda I),
\]
which is equivalent to
\[
Jac(II_{\nu^\perp},I) = 0.
\]
But this equation is the differential equation of
$\nu^\perp$-principal curvature lines (\ref{pcl2}).

Theorem 3.4 of \cite{RS}, lemma 2.1 and theorem 2.1 of \cite{M2}
and above results are put together in the next theorem.

\begin{teo} \label{teo5}
Let $\alpha:\mathbf{M}^2\rightarrow \mathbf R ^4$ be an immersion
of a smooth oriented surface. The following are equivalent
conditions on $\alpha$:
\begin{enumerate}

\item[\textrm a)] The immersion $\alpha$ has everywhere defined
orthogonal asymptotic lines;

\item[\textrm b)] The normal curvature of $\alpha$ vanishes at
every point;

\item[\textrm c)] The immersion $\alpha$ is $\nu$-umbilical for
some unitary normal vector field $\nu$;

\item[\textrm d)] All points of $\alpha$ are semiumbilic;

\item[\textrm e)] There exist a positive scalar function $\lambda$
and an unitary normal vector field $\nu$ such that the second
fundamental form relative to $\nu$ is given by $II_\nu = \lambda
I$;

\item[\textrm f)] The asymptotic lines coincide with the lines of
axial curvature defined by the large axis of the ellipse of
curvature;

\item[\textrm g)] The asymptotic lines coincide with the
$\nu^\perp$-principal curvature lines, for some unitary normal
vector field $\nu$;

\item[\textrm h)] The quartic differential equation of lines of
axial curvature is the product of the quadratic differential
equations of mean directionally curved lines and asymptotic lines.

\end{enumerate}
\noindent  Furthermore if the above function $\lambda$ is a
nonzero constant and the Gaussian curvature $K \neq \lambda^2$
then $\alpha$ is hyperspherical.

\end{teo}

\section{%4.
Concluding remarks}\label{S:4} One direction of research can be
stated: To give an example of a non-hyperspherical immersion
$\alpha$ of a smooth oriented surface in $\mathbf R ^4$ with
globally defined orthogonal asymptotic lines having an isolated
inflection point.

Other direction of research emerges with the evaluation of the
index of an isolated $\nu$-umbilic point. This is related to the
upper bound 1 for the umbilic index on surfaces immersed in
$\mathbf R ^3$ and the Carath\'eodory conjecture (see \cite{SM}
and references therein). Gutierrez and S\'anchez-Bringas
\cite{GuS} have shown that this bound does not hold for the $\nu$
approach.

\end{document}